\documentclass[a4paper,12pt]{amsart}
\usepackage{amssymb}
\usepackage{ifthen}
\usepackage{graphicx}
\usepackage{float}
\usepackage{caption}
\usepackage{subcaption}
\usepackage{cite}
\usepackage{amsfonts}
\usepackage{amscd}
\usepackage{amsxtra}
\usepackage{color}

\newtheorem{thm}{Theorem}
\newtheorem{cor}{Corollary}
\newtheorem{lem}{Lemma}
\newtheorem{rem}{Remark}

\newtheorem{conj}{Conjecture}
\newtheorem{prob}{Problem}

\theoremstyle{definition}
\newtheorem{defn}{Definition}[section]
\newtheorem{example}{Example}

\newenvironment{pf}[1][]{%
 \vskip 1mm
 \noindent
 \ifthenelse{\equal{#1}{}}%
  {{\slshape Proof. }}%
  {{\slshape #1.} }%
 }%
{\qed\bigskip}

\newcounter{alphabet}
\newcounter{tmp}

\makeatletter
\newcommand{\Ref}[1]{\@ifundefined{r@#1}{}{\setcounter{tmp}{\ref{#1}}\Alph{tmp}}}
\makeatother

\newcommand{\IC}{{\mathbb C}}
\newcommand{\ID}{{\mathbb D}}

\def\be{\begin{equation}}
\def\ee{\end{equation}}

\newcommand{\bee}{\begin{enumerate}}
\newcommand{\eee}{\end{enumerate}}

\newcommand{\blem}{\begin{lem}}
\newcommand{\elem}{\end{lem}}
\newcommand{\bthm}{\begin{thm}}
\newcommand{\ethm}{\end{thm}}
\newcommand{\bcor}{\begin{cor}}
\newcommand{\ecor}{\end{cor}}
\newcommand{\beg}{\begin{example}}
\newcommand{\eeg}{\end{example}}
\newcommand{\begs}{\begin{examples}}
\newcommand{\eegs}{\end{examples}}
\newcommand{\bdefe}{\begin{defn}}
\newcommand{\edefe}{\end{defn}}
\newcommand{\bprob}{\begin{prob}}
\newcommand{\eprob}{\end{prob}}
\newcommand{\bques}{\begin{ques}}
\newcommand{\eques}{\end{ques}}
\newcommand{\bei}{\begin{itemize}}
\newcommand{\eei}{\end{itemize}}
\newcommand{\bcon}{\begin{conj}}
\newcommand{\econ}{\end{conj}}
\newcommand{\bcons}{\begin{conjs}}
\newcommand{\econs}{\end{conjs}}
\newcommand{\bprop}{\begin{propo}}
\newcommand{\eprop}{\end{propo}}
\newcommand{\br}{\begin{rem}}
\newcommand{\er}{\end{rem}}
\newcommand{\brs}{\begin{rems}}
\newcommand{\ers}{\end{rems}}
\newcommand{\bo}{\begin{obser}}
\newcommand{\eo}{\end{obser}}
\newcommand{\bos}{\begin{obsers}}
\newcommand{\eos}{\end{obsers}}
\newcommand{\bpf}{\begin{pf}}
\newcommand{\epf}{\end{pf}}
\newcommand{\ba}{\begin{array}}
\newcommand{\ea}{\end{array}}
\newcommand{\beq}{\begin{eqnarray}}
\newcommand{\beqq}{\begin{eqnarray*}}
\newcommand{\eeq}{\end{eqnarray}}
\newcommand{\eeqq}{\end{eqnarray*}}

\newcommand{\ra}{\rightarrow}

\newcommand{\ds}{\displaystyle}

\newcounter{minutes}
\divide\time by 60
\newcounter{hours}
\multiply\time by 60 \addtocounter{minutes}{-\time}

\begin{document}
\bibliographystyle{amsplain}
\title[Bohr radius for locally univalent harmonic mappings]{Bohr radius for locally univalent  harmonic mappings}

\thanks{
File:~\jobname .tex,
          printed: \number\day-\number\month-\number\year,
          \thehours.\ifnum\theminutes<10{0}\fi\theminutes}


\author{Ilgiz R Kayumov, Saminathan Ponnusamy, Nail Shakirov}

\address{I. R Kayumov, Kazan Federal University, Kremlevskaya 18, 420 008 Kazan, Russia
}
\email{ikayumov@kpfu.ru}

\address{S. Ponnusamy, Stat-Math Unit,
Indian Statistical Institute (ISI), Chennai Centre,
110, Nelson Manickam Road,
Aminjikarai, Chennai, 600 029, India.}
\email{samy@isichennai.res.in, samy@iitm.ac.in}

\address{N. Shakirov, Kazan Federal University, Kremlevskaya 18, 420 008 Kazan, Russia
}




\begin{abstract}
We consider the class of all sense-preserving harmonic mappings $f= h+\overline{g}$ of the unit disk $\ID$, where $h$ and
$g$ are analytic with $g(0)=0$, and determine the Bohr radius if any one of the following conditions holds:
\bee
\item $h$ is bounded in $\ID$.
\item $h$ satisfies the condition ${\rm Re}\, h(z)\leq  1$ in $\mathbb{D}$ with $h(0)>0$.
\item both $h$ and $g$ are bounded in $\ID$.
\item $h$ is bounded and $g'(0)=0$.
\eee
We also consider the problem of determining the Bohr radius when the supremum of the modulus of the dilatation of $f$ in $\ID$
is strictly less than $1$. In addition, we determine the Bohr radius for the space $\mathcal B$
of analytic Bloch functions and the space  ${\mathcal B}_H$ of harmonic Bloch functions. The paper concludes with two conjectures.
\end{abstract}


\maketitle
\pagestyle{myheadings}
\markboth{I. R. Kayumov, S. Ponnusamy and N. Shakirov}{Bohr radius for locally univalent  harmonic mappings}

\section{Introduction and Preliminaries}

We shall investigate Bohr's radius for complex-valued harmonic mappings and locally univalent harmonic mappings defined
on the unit disk $\ID :=\{z\in\IC:\, |z|<1\}$. The Bohr
theorem about power series (after subsequent improvements due to M.~Riesz, I.~Schur and F.~Wiener) states that if $f$ is a bounded analytic function
on $\ID$, with the Taylor expansion $f(z)=\sum_{n=0}^\infty a_n z^n$, then
$$\sum_{n=0}^\infty |a_n|r^n \leq \|f\|_\infty
$$
for $0<r\leq 1/3$ and the constant $1/3$ is sharp.
The best constant $r$ in the above inequality, which is $1/3,$ is called the Bohr radius for the class of all analytic self-maps of the unit disk $\ID$.
The original problem goes back to 1914's. Many mathematicians have contributed toward the understanding of
this problem in several settings. We refer to the recent survey on this topic by Abu-Muhanna et al. \cite{AAPon1} for the importance, background, and several
other recent results and extensions. For certain recent results, see \cite{AliBarSoly,KayPon1, KayPon2}.

A harmonic mapping in $\ID$ is a complex-valued function $f=u+iv$ of $z=x+iy$ in $\ID$, which satisfies the Laplace equation
$\Delta f=4f_{z\,\overline{z}}=0$, where $f_{z}=(1/2)\big( f_x-if_y\big)$ and $f_{\overline{z}}=(1/2)\big(f_x+if_y\big)$
and where $u$ and $v$ are real-valued harmonic functions on $\ID$. It follows that $f$ admits the canonical representation
$f= h+\overline{g}$, where $h$ and $g$ are analytic in $\ID$ with $f(0)=h(0)$.
The Jacobian $J_{f}$ of $f$ is given by $J_{f} =|h'|^2-|g'|^2.$ We say that $f$ is sense-preserving in $\ID$ if $J_{f}(z)>0$ in $\ID$. Consequently,
$f$ is locally univalent and sense-preserving in $\ID$ if and only if
$J_{f}(z)>0$ in $\ID$; or equivalently if $h'\neq 0$ in $\ID$ and the dilatation
$\omega_f=:\omega =g'/h'$ has the property that $|\omega (z)|<1$ in $\ID$ (see \cite{Lewy}).

In order to state the first result about Bohr radius for quasiconformal harmonic mappings, we need to introduce some notation.
For harmonic mappings $f$ in $\mathbb{D}$, we use the following
standard notations:
$$\Lambda_{f}(z)=\max_{0\leq \theta\leq 2\pi}|f_{z}(z)+e^{-2i\theta}f_{\overline{z}}(z)|
=|f_{z}(z)|+|f_{\overline{z}}(z)|
$$
and
$$\lambda_{f}(z)=\min_{0\leq \theta\leq 2\pi}|f_{z}(z)+e^{-2i\theta}f_{\overline{z}}(z)|
=\big | \, |f_{z}(z)|-|f_{\overline{z}}(z)|\, \big |
$$
so that if $f$ is locally univalent and sense-preserving, then
$$J_f=\lambda_{f}\Lambda_{f}=|f_{z}|^2-|f_{\overline{z}}|^2 >0.
$$

A sense-preserving homeomorphism $f$ from the unit disk $\ID$ onto $\Omega'$, contained in the Sobolev class $W_{loc}^{1,2}(\ID)$,
is said to be a {\it $K$-quasiconformal mapping} if, for $z\in\ID$,
$$\frac{|f_{z}|+|f_{\overline{z}}|}{|f_{z}|-|f_{\overline{z}}|}= \frac{1+|\omega_f(z)|}{1-|\omega_f(z)|}\leq K,
 ~\mbox{ i.e., }~|\omega_f(z)|\leq k=\frac{K-1}{K+1},
$$
where $K\geq1$ so that $k\in [0,1)$ (cf. \cite{LV,V}).

\begin{thm} \label{KSP4-th3}
Suppose that $f(z) = h(z)+\overline{g(z)}=\sum_{n=0}^\infty a_n z^n+\overline{\sum_{n=1}^\infty b_n z^n}$ is a sense-preserving
$K$--quasiconformal  harmonic mapping of the disk $\ID$, where $h$ is a bounded function in $\ID$. Then
$$\sum_{n=0}^\infty |a_n|r^n+\sum_{n=1}^\infty |b_n|r^n \leq ||h||_{\infty}~\mbox{ for  }~r \leq  \frac{K+1}{5K+1}.
$$
The constant $(K+1)/(5K+1)$ is sharp.
\end{thm}


\begin{thm}\label{cor-ex1}
Assume the hypothesis of Theorem {\rm \ref{KSP4-th3}}. Then
$$|a_0|^2+\sum_{n=1}^\infty (|a_n|+|b_n|)r^n \leq ||h||_{\infty}~\mbox{ for  }~ r \leq \frac{K+1}{3K+1}.
$$
The constant $(K+1)/(3K+1)$ is sharp.
\end{thm}

We would like to remark that the boundedness condition on $h$ in Theorem \ref{KSP4-th3} can be replaced by half-plane
condition. However, the Bohr radius remains the same in this case too.

\begin{thm} \label{KSP4-th4}
Suppose that $f(z) = h(z)+\overline{g(z)}=\sum_{n=0}^\infty a_n z^n+\overline{\sum_{n=1}^\infty b_n z^n}$ is a sense-preserving
$K$--quasiconformal  harmonic mapping of the disk $\ID$, where  $h$ satisfies the condition ${\rm Re}\, h(z)\leq  1$
in $\mathbb{D}$ and $h(0)=a_0$ is positive. Then
$$a_0+\sum_{n=1}^\infty |a_n|r^n+\sum_{n=1}^\infty |b_n|r^n \leq 1 ~\mbox{ for  }~r \leq  \frac{K+1}{5K+1}.
$$
The constant $(K+1)/(5K+1)$ is sharp.
\end{thm}

The following corollaries are regarded as harmonic analogs of the classical Bohr inequality and will be of independent interest.
These results are obtained by allowing $K\ra \infty$ in the above results.


\bcor
Suppose that $f(z) = h(z)+\overline{g(z)}=\sum_{n=0}^\infty a_n z^n+\overline{\sum_{n=1}^\infty b_n z^n}$ is a
sense-preserving harmonic mapping of the disk $\ID$, where $h$ is a bounded function in $\ID$. Then
\be\label{KSP4-eq1}
|a_0| +\sum_{n=1}^\infty ( |a_n|+|b_n|)r^n \leq ||h||_{\infty} ~\mbox{ for  }~ r \leq \frac{1}{5},
\ee
and the number $1/5$ is sharp. Moreover, either $a_0=0$ or $|a_0|$ in \eqref{KSP4-eq1} is replaced by $|a_0|^2$,
then the constant $1/5$ could be replaced by $1/3$ which is also sharp.
\ecor

\bcor
Suppose that $f(z) = h(z)+\overline{g(z)}=\sum_{n=0}^\infty a_n z^n+\overline{\sum_{n=1}^\infty b_n z^n}$ is a
sense-preserving harmonic mapping of the disk $\ID$, where  $h$ satisfies the condition ${\rm Re}\, h(z)\leq  1$
in $\mathbb{D}$ and $h(0)=a_0$ is positive.  Then
$$a_0 +\sum_{n=1}^\infty ( |a_n|+|b_n|)r^n \leq 1 ~\mbox{ for  }~ r \leq  \frac{1}{5},
$$
and the number $1/5$ is sharp.
\ecor

%

\bthm\label{KSP4-th2}
Suppose that either $f=h+g$ or $f = h+\overline{g}$, where $h(z)=\sum_{n=1}^\infty a_n z^n$ and $g(z)=\sum_{n=1}^\infty b_n z^n$
are bounded analytic functions in $\ID$. Then
$$\sum_{n=1}^\infty (|a_n|+|b_n|)r^n \leq \max\{||h||_{\infty}, ||g||_{\infty}\} ~\mbox{ for  }~r \leq \sqrt{\frac{7}{32}}.
$$
This number $\sqrt{7/32}$ is sharp.
\ethm

As in the symmetric case of analytic functions (see \cite{AliBarSoly,KayPon1,KayPon2}), we have the following analog result for
harmonic functions.

\bthm\label{KSP4-Add1}
Let $p \ge 2$.
Suppose that $f(z) = h(z)+\overline{g(z)}=\sum_{n=0}^\infty a_n z^{pn+1}+\overline{\sum_{n=0}^\infty b_n z^{pn+1}}$ is a harmonic $p$--symmetric function
in $\ID$, where $h$ and $g$ are bounded functions in $\ID$. Then
$$\sum_{n=0}^\infty (|a_n|+|b_n|)r^{pn+1} \leq \max\{||h||_{\infty}, ||g||_{\infty}\} ~\mbox{ for  }~r \leq  \frac{1}{2}.
$$
The number $1/2$ is sharp.
\ethm

The proofs of these results will be given in Section \ref{sec2}. In Section \ref{sec3}, we extend further results for
sense-preserving $K$--quasiconformal harmonic mappings of the disk $\ID$. In Section \ref{sec4}, we consider the problem of
finding the Bohr radius for the space of bounded harmonic Bloch functions. The paper concludes with a couple of conjectures.

\section{The Proofs of Theorems \ref{KSP4-th3}, \ref{cor-ex1}, \ref{KSP4-th4},
 \ref{KSP4-th2}  and \ref{KSP4-Add1}
}\label{sec2}

The following lemma is needed for the proof our first theorem.

\begin{lem} \label{KSP4-lem1}
Suppose that $h(z)=\sum_{n=0}^\infty a_n z^n$ and $g(z)=\sum_{n=0}^\infty b_n z^n$ are two analytic functions
in the unit disk $\ID$ such that $|g'(z)| \leq k |h'(z)|$ in $\ID$ and for some $k\in [0,1]$. Then
$$\sum_{n=1}^\infty |b_n|^2r^n \leq k^2\sum_{n=1}^\infty |a_n|^2r^n ~\mbox{ for }~ |z|=r<1.
$$
\end{lem}
\bpf
We integrate inequality $|g'(z)|^2 \leq k^2|h'(z)|^2$ over the circle $|z|=r$ and get
$$
\sum_{n=1}^\infty n^2|b_n|^2 r^{2(n-1)}  \leq  k^2 \sum_{n=1}^\infty n^2|a_n|^2 r^{2(n-1)} .
$$ We integrate the last inequality with respect to $r^2$ and obtain
$$
\sum_{n=1}^\infty n |b_n|^2 r^{2n} \leq   k^2\sum_{n=1}^\infty n|a_n|^2 r^{2n} .
$$
One more integration (after dividing by $r^2$) gives desired inequality.
\epf

\bpf[Proof of Theorem \ref{KSP4-th3}]
For simplicity, we suppose that $||h||_\infty = 1$. Then $|a_n|\leq 1-|a_0|^2$ for $n\geq 1$.
Let $\omega $ denote the dilatation of $f = h+\overline{g}$ so that $|g'(z)| \leq k|h'(z)|$ in $\ID$, where $k\in [0,1)$
and so,
by Lemma \ref{KSP4-lem1},  it follows that
$$
\sum_{n=1}^\infty |b_n|^2r^n\leq  k^2 \sum_{n=1}^\infty |a_n|^2r^n \leq  k^2 (1-|a_0|^2)^2\frac{r}{1-r}.
$$
Consequently,
$$\sum_{n=1}^\infty |b_n|r^n \leq \sqrt{\sum_{n=1}^\infty |b_n|^2r^n}\sqrt{\sum_{n=1}^\infty r^n}
 \leq k(1-|a_0|^2)\frac{r}{1-r}
$$
so that
\be\label{KSP4-eq2}
\sum_{n=0}^\infty |a_n|r^n+\sum_{n=1}^\infty |b_n|r^n \leq |a_0|+(1-|a_0|^2)(1+k)\frac{r}{1-r}
\ee
which is clearly less than or equal to $1$ for $r\leq  1/(3+2k)$. Substituting $k= (K-1)/(K+1)$ gives the desired result.

To prove the sharpness, consider
$$h(z)=\frac{a-z}{1-\overline{a}z} =a + \sum_{n=1}^\infty a_n z^n, \quad a_n=-(1-|a|^2)(\overline{a})^{n-1} \mbox{ for $n\geq 1$},
$$
and $g(z)=\lambda k h(z)$,  where $|\lambda|=1$,   $a \in \ID$ and $k= (K-1)/(K+1)$.
Then it is a simple exercise to see that
$$\sum_{n=0}^\infty |a_n|r^n+\sum_{n=1}^\infty |b_n|r^n =|a| +(1-|a|^2)(1+k)\sum_{n=1}^\infty |a|^{n-1}r^n
=|a| +(1-|a|^2)(1+k)\frac{r}{1-|a|r}
$$
which is bigger than or equal to $1$ if and only if
$$r\geq \frac{1}{1+k+(2+k)|a|}=\frac{K+1}{2K+(3K+1)|a|}.
$$
This shows that the number $(K+1)/(5K+1)$ cannot be improved, since $|a|$ could be chosen so close to $1^{-}$.
\epf

\bpf[Proof of Theorem \ref{cor-ex1}]
Just adopt the method of proof of  Theorem \ref{KSP4-th3}. The desired
conclusion follows if we consider \eqref{KSP4-eq2} and replace the first term $|a_0|$ in \eqref{KSP4-eq2} by $|a_0|^2$. So we omit the details.
\epf


\bpf[Proof of Theorem \ref{KSP4-th4}] We recall that if $p(z)=\sum_{n=0}^\infty p_{n} z^{n}$ is
analytic in $\ID$ such that ${\rm Re}\, p(z)> 0$ in $\ID$, then $|p_n|\leq 2{\rm Re}\,p_0$
for all $n\geq 1$. Applying this result to $p(z)=1-f(z)$ leads to $|a_{n}|\leq 2(1-a_0)$ for all $n\geq 1$. Thus,
as in the proof of Theorem \ref{KSP4-th3}, we can easily obtain from Lemma \ref{KSP4-lem1} that
$$\sum_{n=1}^\infty |b_n|^2r^n\leq  k^2 \sum_{n=1}^\infty |a_n|^2r^n \leq  4k^2 (1-a_0)^2\frac{r}{1-r}
$$
and
$$\sum_{n=1}^\infty |b_n|r^n \leq 2k(1-a_0)\frac{r}{1-r}
$$
so that
\be\label{KSP4-eq2a}
\sum_{n=0}^\infty |a_n|r^n+\sum_{n=1}^\infty |b_n|r^n \leq a_0+2(1+k)(1-a_0)\frac{r}{1-r}
\ee
which is clearly less than or equal to $1$ for $r \leq 1/(3+2k)$. Again, substituting $k= (K-1)/(K+1)$ gives the desired result.
Moreover, sharpness can be seen by considering functions of the form
$$h(z)=\frac{a-z}{1-az}, ~0<a<1, ~\mbox{ and $\ds g(z)= kh(z)=\frac{K-1}{K+1}h(z)$}.
$$
The proof is complete.
\epf


For the proof of Theorem \ref{KSP4-th2}, without loss of generality, we may assume that
$$\max\{||h||_{\infty}, ||g||_{\infty}\}=1.
$$
Then, $||h||_{\infty}\leq 1$ and it follows from the classical Schwarz inequality that
$$\sum_{n=1}^\infty |a_n|r^n \leq \sqrt{\sum_{n=1}^\infty |a_n|^2}\sqrt{\sum_{n=1}^\infty r^{2n}} \leq \frac{r}{\sqrt{1-r^2}}.
$$
Again, since $||g||_{\infty}\leq 1$, the same inequality is valid for $b_n$. Thus, by combining the resulting
inequality with the last inequality, we find that
$$\sum_{n=1}^\infty (|a_n|+|b_n|)r^n \leq  \frac{2r}{\sqrt{1-r^2}} \leq 1 ~\mbox{ for  }~r \leq
 \frac{1}{\sqrt{5}}.
$$
Although this simple approach gives a good estimate, the number $1/\sqrt{5}$ is not sharp. In order to obtain the sharp estimate
we will use a recent approach of Kayumov and Ponnusamy \cite{KayPon1, KayPon2} which, in particular, settled the problem Ali et al. \cite{AliBarSoly}
on the Bohr radius for odd analytic functions.

\begin{lem}\label{KSP4-lem2}
Suppose $p$ is a natural number and $2r^{2p}<1$.
If $h(z)=\sum_{n=0}^{\infty} a_{pn+1} z^{pn+1}$ is analytic and $|h(z)|\le 1$ for $z\in\ID$, then the following inequalities hold:
\begin{equation}\label{KP2-eq3}
\sum_{n=0}^\infty |a_{pn+1}|r^{pn+1} \leq \left \{ \begin{array}{rl}
\frac{1}{r^{p-1}} (3-2 \sqrt{2}\sqrt{1-r^{2p}}) & \mbox{ for $|a_1| \ge r^p$}\\
2r^{p+1} & \mbox{ for $|a_1| < r^p$.}
\end{array}
\right .
\end{equation}
If $r^p \leq 1/3$, then always
$$ \sum_{n=0}^\infty |a_{pn+1}|r^{pn+1} \leq \max\{2r^{p+1},r\}
$$
holds.
\end{lem}
\bpf
This is a special case of the proof of Theorem 1 in \cite{KayPon1, KayPon2}. So, we omit the details to avoid repetition.
\epf


\bpf[Proof of Theorem \ref{KSP4-th2}]
Since $3-2 \sqrt{2}\sqrt{1-r^2}=1/2$ gives $r=\sqrt{7/32}$ and $2r^2=14/32<1/2$, it follows from Lemma \ref{KSP4-lem2} that
$$\sum_{n=0}^\infty |a_{n+1}|r^{n+1} \leq  \frac{1}{2} ~\mbox{ for  }~r \leq \sqrt{\frac{7}{32}}
$$
and the same inequality is valid for the coefficients $b_n$:
$$\sum_{n=0}^{\infty}|b_{n+1}|r^{n+1} \leq  \frac{1}{2}~\mbox{ for  }~r \leq \sqrt{\frac{7}{32}}.
$$
Summing these two inequalities immediately proves that
$$
\sum_{n=1}^{\infty}(|a_{n}|+|b_{n}|)r^{n} \leq 1~\mbox{ for  }~r \leq \sqrt{\frac{7}{32}}.
$$
To prove that the number $\sqrt{7/32}$ cannot be replaced by a larger one, we let $a=3/{\sqrt{14}}$, $\lambda \in\partial \ID$
and consider
$$h(z)=z\left (\frac{a-z}{1-az}\right )=az -(1-a^2)z^2 \sum_{n=2}^\infty  (az)^{n-2} ~\mbox{ and }~ g(z)= \lambda h(z) .
$$
Then, we find that
$$\sum_{n=1}^\infty (|a_n|+|b_n|)r^n = 2\left ( ar +  \frac{(1-a^2)r^2}{1-ar}\right ) =1 ~\mbox{ for  }~r = \sqrt{\frac{7}{32}}
$$
which shows that the number $\sqrt{7/32}$ cannot be improved.
\epf

\bpf[Proof of Theorem \ref{KSP4-Add1}]
Let $r=1/2$. Since $r^p<1/3$ for $p \ge 2$, it follows from Lemma \ref{KSP4-lem2} and the hypothesis that
$$\sum_{n=0}^\infty |a_{pn+1}|r^{pn+1} \leq \frac{1}{2} ~\mbox{ for  }~r \leq  \frac{1}{2}
$$
and the same inequality is valid for the coefficients $b_n$.
%
As a consequence of adding these two inequalities, we obtain that
$$\sum_{n=0}^{\infty}(|a_{pn+1}|+|b_{pn+1}|)r^{pn+1} \leq 1~\mbox{ for  }~r \leq  \frac{1}{2}.
$$
The function $f(z)=z+\overline{z}$ shows that $1/2$ is sharp.
\epf

\section{Further results on Bohr radius for quasiconformal harmonic mappings}\label{sec3}

\begin{thm} \label{KSP4-th5}
Suppose that $f(z) = h(z)+\overline{g(z)}=\sum_{n=0}^\infty a_n z^n+\overline{\sum_{n=2}^\infty b_n z^n}$ is a
sense-preserving $K$--quasiconformal harmonic mapping of the disk $\ID$, where $h$ is a bounded function in $\ID$. Then
\be\label{KSP4-eq5}
\sum_{n=0}^\infty |a_n|r^n+\sum_{n=2}^\infty |b_n|r^n \leq ||h||_{\infty}~\mbox{ for  }~r \leq r_K,
\ee
where $r_K$ is the  positive root of the equation $M_K(r)=1/2$ and
\be\label{KSP4-eq6}
M_K(r)=\frac{r}{1-r} + \left (\frac{K-1}{K+1}\right ) r^2\sqrt{\frac{1+r^2}{(1-r^2)^3}}\sqrt{\frac{\pi^2}{6}-1}.
\ee
The number $r_K$ cannot be replaced by a number greater than $R=R(K)$, where $R$ is the positive root of the equation
\be\label{KSP4-eq7}
\frac{4R}{1-R}\left (\frac{K}{K+1}\right )+2\left (\frac{K-1}{K+1}\right )\log(1 - R)=1.
\ee
\end{thm}
\bpf
Without loss of generality we may assume that $||h||_{\infty} \leq 1$. Because
$f=h+\overline{g}$ is sense-preserving and $K$--quasiconformal harmonic mapping with $g'(0)=0$,  Schwarz's lemma gives that
$\omega_f=:\omega =g'/h'$ is analytic in $\ID$ and $|\omega (z)|\leq k|z|$ in $\ID$, where $k=(K-1)/(K+1)$.  Thus, we have
$$|g'(z)|^2 = |\omega(z) h'(z)|^2\leq k^2|zh'(z)|^2.
$$
We integrate this inequality over the circle $|z|=r$ and obtain
$$\sum_{n=2}^\infty n^2|b_n|^2 r^{2(n-1)} \leq k^2r^2 \sum_{n=1}^\infty n^2|a_n|^2 r^{2(n-1)}
$$
and, because $|a_n|\leq 1-|a_0|^2$ for $n\geq 1$, we have
$$\sum_{n=2}^\infty n^2|b_n|^2 r^{2n} \leq k^2 (1-|a_0|^2)^2r^4\sum_{n=1}^\infty n^2 r^{2(n-1)}=k^2 (1-|a_0|^2)^2\frac{r^4(1+r^2)}{(1-r^2)^3}.
$$
Consequently, using the classical Schwarz  inequality, we deduce that
\begin{equation} \label{eq1}
\sum_{n=2}^\infty |b_n| r^n \leq \sqrt{\sum_{n=2}^\infty n^2|b_n|^2 r^{2n}}\sqrt{\sum_{n=2}^\infty \frac{1}{n^2}}
\leq kr^2(1-|a_0|^2)\sqrt{\frac{1+r^2}{(1-r^2)^3}}\sqrt{\frac{\pi^2}{6}-1}.
\end{equation}
Therefore, from (\ref{eq1}) it follows that
\beqq
S&=& |a_0|+|a_1|r+\sum_{n=2}^\infty (|a_n|+|b_n|)r^n \\
&\leq&
|a_0|+(1-|a_0|^2)\left [\frac{r}{1-r} + kr^2\sqrt{\frac{1+r^2}{(1-r^2)^3}}\sqrt{\frac{\pi^2}{6}-1}\right ]
\eeqq
which is less than or equal to $1$ provided $M_K(r)=1/2$ holds, where $M_K(r)$ is defined by \eqref{KSP4-eq6}.

Finally, we let $a \in [0,1)$ and consider the functions
$$h(z) = \frac{a-z}{1-az} ~\mbox{ and }~ g'(z)=kzh'(z).
$$
From here we find that
$$ |a_n|=a^{n-1}(1-a^2)~\mbox{ and }~ |b_n|=k\left (\frac{n-1}{n}\right ) a^{n-2}(1-a^2), \quad n \ge 2,
$$
so that
\beqq
\sum_{n=0}^\infty |a_n|r^n+\sum_{n=2}^\infty |b_n|r^n
&=& a+(1-a^2)\sum_{n=1}^\infty a^{n-1}r^n+k(1-a^2)\sum_{n=2}^\infty  \frac{n-1}{n} a^{n-2}r^n\\
&=& a+(1-a^2)\frac{r}{1-ar}+k(1-a^2)\frac{a r + (1 -a r) \log(1 - a r)}{a^2 (1 - a r)}\\
&=& a+(1-a^2)\left [\frac{(a+k)r}{a(1-ar)}+k\frac{\log(1-a r)}{a^2}\right ],
\eeqq
where $k=(K-1)/(K+1)$. Simple analysis shows that the last expression is less than or equal to  $1$ for all $a \in [0,1)$ only
in the case when $r \leq R=R(K)$ which is the positive root of the equation \eqref{KSP4-eq7}.
\epf


Allowing $K\rightarrow \infty$ in Theorem \ref{KSP4-th5} shows that the root $r_\infty$ of the limiting case $M_\infty (r)=1/2$, i.e.,
$$
\frac{r}{1-r} + r^2\sqrt{\frac{1+r^2}{(1-r^2)^3}}\sqrt{\frac{\pi^2}{6}-1}=\frac{1}{2},
$$
gives the value $0.2942...$. We may now formulate this discussion as follows.

\begin{cor} \label{KSP4-cor2}
Suppose that $f(z) = h(z)+\overline{g(z)}=\sum_{n=0}^\infty a_n z^n+\overline{\sum_{n=2}^\infty b_n z^n}$ is a
sense-preserving harmonic mapping of the disk $\ID$, where $h$ is a bounded function in $\ID$. Then
\be\label{KSP4-eq3a}
\sum_{n=0}^\infty |a_n|r^n+\sum_{n=2}^\infty |b_n|r^n \leq ||h||_{\infty}~\mbox{ for  }~r \leq 0.2942... .
\ee
The number $0.2942...$ cannot be replaced by a number greater than $R=0.299825...$, where $R$ is the positive root of the equation
$$\frac{4R}{1-R}+2\log(1 - R)=1.
$$
\end{cor}

Further remarks would be useful.


\br
{\rm
Also, it is worth pointing out that if the first term $|a_0|$ in \eqref{KSP4-eq5} is replaced by $|a_0|^2$, then the conclusion of
Corollary \ref{KSP4-cor2} takes the form
$$|a_0|^2+ \sum_{n=1}^\infty |a_n|r^n+\sum_{n=2}^\infty |b_n|r^n \leq  ||h||_{\infty}~\mbox{ for  }~r \leq r_K
$$
where $r_K$ is the  positive root of the equation $M_K (r)=1$ and $M_K (r)$ is defined by \eqref{KSP4-eq6}.
The number $r_K$ cannot be replaced by a number greater than $R=R(K)$, where $R$ is the positive root of the equation
$$
\frac{2R}{1-R}\left (\frac{K}{K+1}\right )+\left (\frac{K-1}{K+1}\right )\log(1 - R)=1.
$$

Again, the case $K\rightarrow \infty$ (i.e. when the dilatation has the property that $|\omega (z)|<1$ in $\ID$) needs a special mention, since
the corresponding value of $r_\infty$ is $0.435668...$, where the number $0.435668...$ is the root of the equation $M_\infty (r)=1$, i.e.,
$$\frac{r}{1-r} + r^2\sqrt{\frac{1+r^2}{(1-r^2)^3}}\sqrt{\frac{\pi^2}{6}-1}=1.
$$
Furthermore, the number $0.435668...$ cannot be replaced by a number greater than $R(\infty)=0.44182...$, where $R(\infty)$ is the positive root of the equation
$$\frac{2R}{1-R}+\log(1 - R)=1.
$$
}
\er
\begin{thm}\label{KSP4-th8New}
Suppose that $f(z) = h(z)+\overline{g(z)}=\sum_{n=0}^\infty a_n z^n+\overline{\sum_{n=0}^\infty b_n z^n}$ is a
locally univalent $K$--quasiconformal harmonic mapping of the disk $\ID$, where $h'$ is a bounded function in $\ID$. Then
$$
\frac{K+1}{2K}\sum_{n=1}^\infty (|a_n|+|b_n|) nr^{n-1} \leq ||h'||_{\infty} ~\mbox{ for  }~ r \leq 1/3
$$
and the number $1/3$ is sharp.
\end{thm}


The proof of Theorem \ref{KSP4-th8New} easily follows easily if we use the above method and the
classical proof of Bohr's $1/3$-Theorem. The case $K\rightarrow \infty$ gives the corresponding result for
sense-preserving harmonic mappings of the disk $\ID$.

\section{Bohr radius for harmonic Bloch functions}\label{sec4}

A harmonic function $f$ is called a {\it harmonic Bloch function} if and only if
$$\sup_{z\in\mathbb{D}}(1-|z|^{2})\Lambda_{f}(z)<+\infty,
$$
where $\Lambda_{f}(z)=|f_{z}(z)|+|f_{\overline{z}}(z)|.$
The space of all harmonic Bloch functions, denoted by the symbol $\mathcal{B}_H$,
forms a complex Banach space  with the norm $\|\cdot\|$ given by (see \cite{Co})
$$
\|f\|_{\mathcal{B}_H}=|f(0)|+\sup_{z\in\mathbb{D}}(1-|z|^{2})\Lambda_{f}(z),
$$
where $\Lambda_{f}(z) =|f_{z}(z)|+|f_{\overline{z}}(z)|$.  This is referred to as the harmonic Bloch norm and the elements
of the harmonic Bloch space are called harmonic Bloch functions. Recently, the space  $\mathcal{B}_H$ together with its various generalizations
have been studied extensively. See for example, see \cite{CPW-BAMS11,CPW-BMMSS11,CPW-CAOT12}.

Clearly, this definition coincides with the classical Bloch space $\mathcal{B}$ when $f$ is analytic in $\ID$. We refer to the basic paper on
this topic by Anderson et al. \cite{ACP74} and the book of Pommerenke \cite{PommBBCM-92}.

\begin{thm} \label{KSP4-th6}
Let $f\in \mathcal{B}$ and $\|f\|_{\mathcal{B}} \leq 1$. Then
$$ \sum_{n=0}^\infty |a_n| r^n  \leq 1 ~\mbox{ for }~r \leq R=0.55356... ,
$$
where $R$ is the positive solution to the equation
$$1-R+R\log (1-R)=0.
$$
The number $R$ cannot be replaced by a number greater than $0.624162...$.
\end{thm}
\bpf
Let $f(z)=\sum_{n=0}^\infty a_n z^n$ be analytic and
$$||f||_{\mathcal{B}}=|f(0)|+\sup_{z \in \ID}  (1-|z|^2)|f'(z)|\leq 1.
$$
Thus, we have
$$|f'(z)|^2 \leq \frac{(1-|a_0|)^2}{(1-|z|^2)^2}.
$$
We integrate this inequality over the circle $|z|=r$ and obtain
$$\sum_{n=1}^\infty n^2|a_n|^2r^{2n} \leq r^2\frac{(1-|a_0|)^2}{(1-r^2)^2}
$$
so that
$$\sum_{n=1}^\infty n^2|a_n|^2r^{n-1} \leq \frac{(1-|a_0|)^2}{(1-r)^2} ~\mbox{ for $r<1$.}
$$
We may now integrate this with respect to $r$ (with limit from $0$ to $r$) and obtain
$$\sum_{n=1}^\infty n|a_n|^2r^{n-1} \leq \frac{(1-|a_0|)^2}{1-r} ~\mbox{ for $r<1$},
$$
which by integration with respect to $r$ gives
$$\sum_{n=1}^\infty |a_n|^2r^{n} \leq (1-|a_0|)^2\log\frac{1}{1-r} ~\mbox{ for $r<1$.}
$$
Consequently,
$$\sum_{n=1}^\infty |a_n| r^n
\leq  \sqrt{\sum_{n=1}^\infty |a_n|^2 r^{n}}\sqrt{\sum_{n=1}^\infty r^n} \leq (1-|a_0|)\sqrt{\log\frac{1}{1-r}}\sqrt{\frac{r}{1-r}}.
$$
It means that
$$\sum_{n=0}^\infty |a_n| r^n \leq |a_0|+(1-|a_0|)\sqrt{\log\frac{1}{1-r}}\sqrt{\frac{r}{1-r}}
 \leq 1~\mbox{ for }~r \leq R,
$$
where $R=0.55356...$ is the positive solution to the equation
$$\log\frac{1}{1-R}=\frac{1-R}{R}.
$$

Now we consider (see \cite{Avkhkay})
$$f(z)=\frac{3\sqrt{3}}{4}\left(\left(\frac{z-a}{1-a z}\right )^2-a^2\right)=\frac{3\sqrt{3}}{4}\left( \sum_{n=1}^\infty a_n z^n \right),
$$
where $a\in (0,1)$, $a_1= -2a (1 - a^2)$ and
$$ a_n=(1-a^2)[n(1-a^2)-(1+a^2)](a)^{n-2} ~\mbox{ for $n\geq 2$}.
$$
Since
$$f'(z)=\frac{3\sqrt{3}}{2}\left(\frac{z-a}{1-a z}\right)\frac{1-a^2}{(1-a z)^2},
$$
it is easy to check that $||f||_{\mathcal{B}}= 1$. For this function, we observe that the coefficients $a_n$ for $n\geq 2$
are all positive, whenever $a \in (0,1/\sqrt{3})$. Furthermore,
$$\sum_{n=0}^\infty |a_n| r^n=\frac{3\sqrt{3}}{4}\left(\left(\frac{r-a}{1-a r}\right )^2-a^2 +4a(1-a^2)r \right).
$$
Now we  suppose that $r$ is a function of $a$. We want to find the minimal $r$ for which
$$\frac{3\sqrt{3}}{4}\left(\left(\frac{r-a}{1-ar}\right )^2-a^2+4a(1-a^2)r \right)=1
$$
which may be rewritten as
\begin{equation} \label{Extra1}
(r-a)^2+(1-ra)^2(-a^2 + 4a(1 - a^2)r) =\frac{4}{3\sqrt{3}}(1-ar)^2.
\end{equation}
Now we differentiate $r$ in the variable $a$ and then we set $r'(a)=0$. We arrive at the cubic equation
\begin{equation} \label{Extra2}
18 r + 8 \sqrt{3} r - 54 a^2 r - 144 a r^2 - 8 \sqrt{3} a r^2 +
 252 a^3 r^2 + 108 a^2 r^3 - 180 a^4 r^3=0.
\end{equation}
The algebraic system of the equations (\ref{Extra1}) and (\ref{Extra2}) can be easily solved, for example, by Mathematica 10.
Consequently, we get $a=0.3775$ and then we obtain $r=0.624162$ such that
$\sum_{n=0}^\infty |a_n| r^n=1.
$
%
The proof is complete. \epf

\begin{thm}  \label{KSP4-th7}
Suppose that $f=h+\overline{g}$ is harmonic in $\ID$, $g(0)=0$ and $||f||_{\mathcal{B}_H}\leq 1$,
where
$$||f||_{\mathcal{B}_H}=|f(0)|+\sup_{z \in \ID}  (1-|z|^2)(|h'(z)|+|g'(z)|).
$$
Then
$$ |a_0|+\sum_{n=1}^\infty \sqrt{|a_n|^2+|b_n|^2} r^n  \leq 1~\mbox{ for }~ r \leq R=0.55356.
$$
This number $0.55356$ cannot be replaced by a number greater than $0.624162...$.
\end{thm}
\bpf
By assumption, we have
$$ |h'(z)|^2+|g'(z)|^2 \leq \frac{(1-|a_0|)^2}{(1-|z|^2)^2}.
$$ We integrate this inequality over the circle $|z|=r$ and obtain
$$
\sum_{n=1}^\infty n^2(|a_n|^2+|b_n|^2)r^{2n} \leq r^2\frac{(1-|a_0|)^2}{(1-r^2)^2}.
$$
Now the remaining part of the proof is identical to that of Theorem \ref{KSP4-th6}. Thus, the proof is complete.
\epf

As remarked earlier if we replace the first term $|a_0|$ in the conclusion of the last two theorems by $|a_0|^2$, then the Bohr radius
obviously can be stated in an improved form.

%
%
%

\section{Conclusion}
We conclude the paper with the following conjectures.

\bcon\label{conj1}
Suppose that $f = h+\overline{g}=\sum_{n=1}^\infty a_n z^n+\overline{\sum_{n=1}^\infty b_n z^n}$ is a harmonic mapping of
the disk $\ID$, where $h$ is bounded in $\ID$. If $|g'(z)| \leq \frac{K-1}{K+1}|h'(z)|$, then
$$ \sum_{n=1}^\infty (|a_n|+|b_n|)r^n \leq ||h||_{\infty} ~\mbox{ for }~  r \leq \frac{1}{4} \sqrt{\frac{7}{2} - \frac{1}{2K^2} + \frac{5}{K}}.
$$
This constant is sharp for all  $K\geq 1$.
\econ

The conjectured extremal function has the form
$$f(z)=h(z)+\frac{K-1}{K+1}\overline{h(z)}, \quad h(z)=z\left(\frac{z-a}{1-\overline{a}z}\right),
$$
with a suitable $a \in \mathbb{D}$.

If we replace the condition ``$|g'(z)| \leq \frac{K-1}{K+1}|h'(z)|$'' by ``$f$ is $K$-quasiconformal", then the Bohr radius will be greater than
the number mentioned in Conjecture \ref{conj1}, because the conjectured extremal function is not locally univalent in the unit disk.

The proof of Theorem \ref{KSP4-Add1} leads to another problem which we state it now as a conjecture.

\bcon
Let $p \ge 2$.
Suppose that $f(z) = h(z)+\overline{g(z)}=\sum_{n=0}^\infty a_n z^{pn+1}+\overline{\sum_{n=0}^\infty b_n z^{pn+1}}$ is a harmonic
$p$--symmetric and sense-preserving mapping in $\ID$, where $h$ and $g$ are bounded functions in $\ID$. Then
$$\sum_{n=0}^\infty (|a_n|+|b_n|)r^{pn+1} \leq
\max\{||h||_{\infty},||g||_{\infty}\} ~\mbox{ for  }~r \leq \frac{1}{2}.
$$
The constant $1/2$ is sharp.
\econ

Also it would be interesting to obtain an analog of the Conjecture 2 for locally $K$-quasiconformal mappings.

\subsection*{Acknowledgements}
The research of the first author was supported by Russian foundation for basic research, Proj. 17-01-00282,
and the research of the second author was supported
by the project RUS/RFBR/P-163 under Department of Science \& Technology (India).
The second author is currently on leave from the IIT Madras.

\end{document}